\documentclass[twoside,leqno, twocolumn]{article}

\usepackage[letterpaper]{geometry}
\usepackage{amssymb, amsmath, accents}
\usepackage{graphicx, caption, subcaption, tabularx}
\usepackage[url=false]{biblatex}
\usepackage{ltexpprt}
\usepackage{hyperref, cleveref}

\newcommand{\abs}[1]{\ensuremath{\left|#1\right|}}
\DeclareMathOperator*{\argmin}{\arg\min}
\DeclareMathOperator{\card}{\mathrm{card}}
\newcommand{\bbR}{\ensuremath{\mathbb{R}}}

\newcommand{\frS}{\ensuremath{\mathfrak{S}}}
\newcommand\ubar[1]{\underaccent{\bar}{#1}}

\newtheorem{assumption}{Assumption}[section]

\renewbibmacro*{url+urldate}{}

\begin{document}

\newcommand\relatedversion{}

\title{\Large An Efficient Local Optimizer-Tracking Solver for Differential-Algebriac Equations with Optimization Criteria\relatedversion}
% which affiliation should I put here?
\author{Alexander Fleming\thanks{RWTH Aachen University. }
\and Jens Deussen$^*$
\and Uwe Naumann$^*$}

\date{}

\maketitle

% Copyright Statement
% When submitting your final paper to a SIAM proceedings, it is requested that you include
% the appropriate copyright in the footer of the paper.  The copyright added should be
% consistent with the copyright selected on the copyright form submitted with the paper.
% Please note that "20XX" should be changed to the year of the meeting.

% Default Copyright Statement
% \fancyfoot[R]{\scriptsize{Copyright \textcopyright\ 20XX by SIAM\\
% Unauthorized reproduction of this article is prohibited}}

% Depending on which copyright you agree to when you sign the copyright form, the copyright
% can be changed to one of the following after commenting out the default copyright statement
% above.

\fancyfoot[R]{\scriptsize{Copyright \textcopyright\ 2024\\
Copyright for this paper is retained by authors}}

%\fancyfoot[R]{\scriptsize{Copyright \textcopyright\ 20XX\\
%Copyright retained by principal author's organization}}

%\pagenumbering{arabic}
%\setcounter{page}{1}%Leave this line commented out.

% send help
\begin{abstract} \small\baselineskip=9pt 
	A sequential solver for differential-algebraic equations with embedded optimization criteria (DAEOs) was developed to take advantage of the theoretical work done by Deussen et al. in   \cite{deussenNumericalSimulationDifferentialalgebraic2023}. Solvers of this type separate the optimization problem from the differential equation and solve each individually. The new solver relies on the reduction of a DAEO to a sequence of differential inclusions separated by \emph{jump events}. These jump events occur when the global solution to the optimization problem jumps to a new value. Without explicit treatment, these events will reduce the order of convergence of the integration step by one   \cite{mannshardtOnestepMethodsAny1978}.
	The solver implements a "local optimizer tracking" procedure to detect and correct these jump events. Local optimizer tracking is much less expensive than running a deterministic global optimizer at every time step. 
	% "solely on an integrator" is weak
	This preserves the order of convergence of the integrator component without sacrificing performance to perform deterministic global optimization at every time step. The newly developed solver produces correct solutions to DAEOs and runs much faster than sequential DAEO solvers that rely only on global optimization.
\end{abstract}

\section{Introduction}
Differential-algebraic equations with embedded optimization criteria (DAEO) describe dynamic systems where at least one of the algebraic variables is determined by the solution to an optimization problem. Phase-separation problems subject to thermodynamic equilibrium constraints (for example, vapor-liquid separation) are common examples of DAEOs, where the equilibrium between phases occurs at the minimum of the Gibbs free energy   \cite{bakerGibbsEnergyAnalysis1982}.A few specific industrial applications of these models are distillation processes, flash drums, heat exchangers, and bio-chemical reaction processes   \cite{raghunathanMPECFormulationDynamic2004, plochMultiscaleDynamicModeling2019}. In this specific case, the optimality constraint on the Gibbs free energy leads to a system of differential-algebraic equations with an embedded non-linear programming problem, which is equivalent to a DAEO   \cite{bieglerNonlinearProgrammingConcepts2010}.
% maybe go grab the U-V flash drum example? or something else from the references.

Solving these types of problems frequently begins by replacing the embedded optimization problem with the Karush-Kuhn-Tucker (KKT) conditions   \cite{bieglerNonlinearProgrammingConcepts2010}. The KKT conditions state that the solution to the following problem

\begin{equation*}
	\begin{gathered}
		\mathrm{minimize}\,f(x) \\
		\mathrm{subject\ to} \\
		g_i(x) = 0\\
		h_j(x) \leq 0
	\end{gathered}
\end{equation*}
is equivalent to finding all saddle points of
\begin{equation}
	\label{eq:kkt-lagrangian}
	\mathcal{L}(x, \mu, \lambda) = f(x) + \mu_jh_j(x) + \lambda_ig_i(x)
\end{equation}

In problems without constraints, the KKT conditions can be reduced to simpler first order optimality conditions.

The KKT conditions guarantee first-order local optimality, but cannot provide information about global optimality. A good initial guess for solving \eqref{eq:kkt-lagrangian} can, however, arrive at the global optimum. Therefore, using this method can only provide an exact solution when the embedded optimization problem is a convex non-linear programming problem, or if the initially computed global minimizer is and remains the global minimizer. If this is not the case, the evolution of the system dynamics will generate multiple feasible solutions to the problem rewritten using the KKT conditions \cite{plochDirectSingleShooting2022}.

Methods for solving the reformulated problem are divided into two classes: sequential and simultaneous. Simultaneous methods solve the reformulated problem by combining the entire set of equations and variables into one nonlinear programming problem, which can be solved using suitable tools, such as IPOPT   \cite{bieglerNonlinearProgrammingConcepts2010, wachterImplementationInteriorpointFilter2006}. Sequential methods separate the differential-algebraic equation from the optimization problem, and require some way to solve both a DAE and an optimization problem. Sequential approaches have been used to solve DAEOs that describe U-V flash processes in   \cite{ritschelAlgorithmGradientbasedDynamic2018}, and flux-balance problems in   \cite{plochDirectSingleShooting2022}. However, these single-shooting algorithms for solving non-smooth DAEs fail to account for the fact that the KKT conditions do not guarantee global optimality.

This work develops a sequential solver for DAEOs and builds upon the theoretical work done in   \cite{deussenNumericalSimulationDifferentialalgebraic2023}. Rather than applying the KKT conditions and successively identifying the optimizer, the solver finds \textit{all} local optimizers in a given search domain via a branch-and-act algorithm   \cite{bongartzMAiNGOMcCormickbasedAlgorithm, rallGlobalOptimizationUsing1985}. 

Given certain assumed properties of the DAEO, the movement of these local optimizers can be tracked through a time step by appending the first-order optimality conditions (gradient is zero) to the equations solved by the integrator component of the solver. Solving these equations is significantly less expensive than running a global optimizer at every time step, and it is possible to correct for possible "jump events", which occur when the identity of the global optimizer changes \textit{during} a time step. If left unaccounted for, these jump events will reduce the order of convergence of the solution to the DAE to 1 \cite{deussenNumericalSimulationDifferentialalgebraic2023, mannshardtOnestepMethodsAny1978}.

The solver developed in this work performs both local optimizer tracking and jump event correction. The numerical experiments performed in \Cref{section:numerical-experiments} demonstrate the correctness of local optimizer tracking and the performance improvement over a sequential DAEO solver without local optimizer tracking.
% can't back this up, not without a comparison to an NLP solver.
% These processes are much more computationally efficient than solving a (potentially) large non-linear programming problem at each time step.

\section{Problem Description.}
A simple DAEO is written as an initial-value problem like so:
\begin{equation} \label{eq:daeo-ivp}
\begin{aligned}
	x(t_0) &= x_0\\
	\dot{x}(x, y) &= f(x(t), y(t))\\
	\left\{y^k\right\}(x, t) &= \argmin_{y}h(x(t), y)
\end{aligned}
\end{equation}
where 
\begin{itemize}
	\item $f:\bbR^{n_x}\times\bbR^{n_y}\mapsto\bbR^{n_x}$ describes the differential behavior of the system
	\item $h:\bbR^{n_x}\times\bbR^{n_y}\mapsto\bbR$ is an objective function to be minimized
	\item $\left\{y^k(t)\right\}$ denotes the set of minimizers of $h(x, y)$ at time $t>t_0$
\end{itemize}
The notation $\dot{x}$ denotes $\partial x/\partial t$, $\partial_x f$ denotes $\partial f/\partial x$, and $d_x f$ denotes $df/dx$, in the cases where the difference between the partial and total derivatives is relevant. Second derivatives are written similarly, with $\ddot{x}$ for $\partial^2_{tt} x$ and $\partial^2_{yx} f$ for $\partial f/\partial x\partial y$.

At any time $t \geq t_0$, there will be some $y^\star(t)\in\left\{y_k(t)\right\}$ such that
\begin{equation}
	h(x(t), y^\star(t))\leq h(x(t),y^i(t))\, \forall\,y^i(t)\in\left\{y^k(t)\right\} 
\end{equation}

This implies that the initial value problem posed in \eqref{eq:daeo-ivp} can be divided into a sequence of initial-value problems if $y^\star(t)$ has a sufficiently well-behaved dependency on $x$ and $t$  \cite{deussenNumericalSimulationDifferentialalgebraic2023}.

\begin{assumption}
	\label{assume:events-exist}
	There is only a finite set of events $t_0 < \tau_1 < \tau_2 < \ldots$ where more than one possible global optimizer exists:
	\begin{equation*}
		\card\left(\argmin_{y\in\left\{y^k(\tau_j)\right\}}h(x(\tau_j), y)\right) > 1
	\end{equation*}

	\begin{equation*}
		\argmin_{y\in\left\{y^k(\tau_j)\right\}} h(x(\tau_j), y) = \left\{y^1, y^2\right\}
	\end{equation*}
\end{assumption}
this assumption limits the behavior at each $\tau_j$ to the following two possibilities:
\begin{enumerate}
	\item Some number of local optimizers $y$ are each global optimizers \textit{only} at $\tau_j$, and $y^1$ is the global optimizer both before and after $\tau_j$:
	\begin{equation*}
		\begin{aligned}
			\lim_{t\to \tau_j^+} \argmin_{y} h(x(t), y) &= y^1\\
			\lim_{t\to \tau_j^-} \argmin_{y} h(x(t), y) &= y^1
		\end{aligned}
	\end{equation*}
	\item $y^1$ is the global optimizer \textit{before} $\tau_j$, and $y^2$ is the global optimizer \textit{after} $\tau_j$.
	\begin{equation*}
		\begin{aligned}
			\lim_{t\to \tau_j^+} \argmin_{y} h(x(t), y) &= y^1\\
			\lim_{t\to \tau_j^-} \argmin_{y} h(x(t), y) &= y^2
		\end{aligned}
	\end{equation*}
\end{enumerate}
In both cases, integrating the differential part of \eqref{eq:daeo-ivp} from $t_1 < \tau_1$ to $t_2 > \tau_1$ yields
\begin{equation}
	\begin{aligned}
		x(t_2) &= x(t_1) + \int_{t_1}^{t_2}f(x, y^\star)\,dt\\
		&=x(t_1) + \int_{t_1}^{\tau_1}f(x, y^1)\,dt + \int_{\tau_1}^{t_2}f(x, y^2)\,dt
	\end{aligned}
\end{equation}
Therefore, the choice of global optimizer $y^\star$ exactly at $\tau_j$ is inconsequential. This renders the first case in \Cref{assume:events-exist} trivially easy to evaluate, since there is no change in the identity of $y^\star$ after the event.

The second case, however, is relevant for the correct solution of DAEOs. Between each $\tau_j$, the DAEO \eqref{eq:daeo-ivp} can be treated as an IVP, with the extra constraint that $\partial_{y} h(x(t), y^\star(t)) = 0$   \cite{deussenSubdomainSeparabilityGlobal2023}.

\begin{assumption} 
	\label{assume:events-transversal}
	(Events Are Transversal)
	The identity of the global optimizer $y^\star$ at the time step $t^{n+1}$ \textit{after} an event is always different from the identity of the global optimizer $y^\star$ at the time step $t^n$ \textit{before} an event. 
\end{assumption}

This transversality assumption eliminates the case that the global optimizer $y^\star$ may switch to a new value and return \textit{inside} a time step. This behavior would occur on a scale that is too small for a numerical solver to capture. If this behavior were relevant to the solution, simply adjusting the time step size such that \Cref{assume:events-transversal} would generate a correct solution.

\begin{assumption}
	\label{assume:global-does-not-emerge}
	A new global optimizer of $h(x, y)$ does not emerge during a time step. Events can only occur between already-existing local optimizers.
\end{assumption}

This assumption is critical for the local optimizer tracking strategy. Without it, the global optimizer would necessarily run every time step to mitigate the possibility of a newly-emerged global optimizer.

% I would like to comment on the reality / feasibility of this assumption. Maybe comment on phase transition? My metallurgy / chemistry knowledge is limited and I'm not sure where to cite. 

To solve this sequence of IVPs generated by the jump events, the solver must 
\begin{enumerate}
	\item Generate the set $\left\{y^k(t)\right\}$ at $t=t_0$ and at user-specified times $t>t_0$ to account for the possible emergence or disappearance of local optima.
	\item Detect when the identity of $y^\star$ changes and correct the time step in which it happens (the second case, above).
	\item Solve the sequence of IVPs generated by jump events for $t>t_0$. The global optimizer is not necessary, unless optimizers of $h(x, y)$ can vanish or emerge as $t$ increases.
\end{enumerate}

\section{Deterministic Global Optimization with Interval Arithmetic.}
A sequential DAEO solver must first locate the global optimizer $y^\star$ of $h(x, y)$. This can be accomplished using a branch-and-act algorithm to locate each point $y^k$ inside an initial search domain that satisfies certain optimality conditions.

\begin{assumption} 
	\label{assume:twice-lipschitz}
	Both $f$ and $h$ are twice Lipschitz continuous with respect to $(x(t), y(t))$.
\end{assumption}
Optimality conditions for $h(x(t), y)$ follow from \Cref{assume:twice-lipschitz}
\begin{equation}
	\label{eq:optimality-conditions}
	\begin{aligned}
		0 &= \partial_yh(x, y)\\
		0 &\prec\partial^2_{y}h(x, y)
	\end{aligned}
\end{equation}
Using a suitable algorithmic differentiation (AD) tool and an interval arithmetic library allows for the comparatively inexpensive verification of the optimality conditions \cite{rallGlobalOptimizationUsing1985}. 

An interval $[y_i] = \left[\ubar y_i, \bar y_i\right]$ bounds all possible values for each $y_i\in\bbR$. 
\begin{Definition}
	\label{def:interval-arith}
	\textbf{(Fundamental Theorem of Interval Arithmetic)}
	
	An interval evaluation of a function $\left[\ubar{y},\bar{y}\right] = f(\left[\ubar x, \bar x\right])$ must yield an interval that contains all pointwise evaluations $f(x)\,\forall\,x\in\left[\ubar x, \bar x\right]$.
\end{Definition}

A consequence of \ref{def:interval-arith} is that a correct implementation of interval arithmetic need not provide the smallest possible result $\left[\ubar y, \bar y\right]$  \cite{hickeyIntervalArithmeticPrinciples2001}. However, the result $[\ubar y, \bar y] = f(\left[\ubar x, \bar x\right])$ \textit{must} contain all possible results of the interval evaluation. One such implementation is the Boost Interval library  \cite{melquiondBoostIntervalLibrary2022}.

AD applies to the chain rule to compute derivatives of arbitrary functions by decomposing them into a sequence of fundamental operations with known derivatives \cite{griewankEvaluatingDerivativesPrinciples2008}. Over the real numbers, tangent-mode AD of some function $g(x)$ would effectively compute both $g$ and its first partial derivative $\partial_x(g)\dot{x}$ from an input $x$ and a perturbation $x^{\left(1\right)}$. Over intervals, however, tangent-mode AD of $g(\left[x\right])$ would compute interval enclosures for both $g(\left[x\right])$ and $\partial_xg(\left[x\right])x^{\left(1\right)}$. 

Similarly to tangent mode, an adjoint-mode AD evaluation of $g(x)$ computes $x_{\left(1\right)}\partial_xg(x)$ from an input $x$ and a sensitivity $x_{\left(1\right)}$. Over intervals, adjoint-mode AD computes an interval enclosure for each of $g(\left[x\right])$ and $x_{\left(1\right)}\partial_xg(\left[x\right])$.

If the gradient $\partial_{y_i}h(x, y)$ must be in an interval $\left[\partial_{y_i}h(x, \left[y\right])\right]$, and that interval does not contain $0$, then the input interval $\left[y\right]$ contains no points that satisfy the optimality conditions \eqref{eq:optimality-conditions}. These interval gradients are computed by substituting appropriately-typed intervals directly into the gradient driver for the function   \cite{griewankEvaluatingDerivativesPrinciples2008, rallGlobalOptimizationUsing1985}. 

Similarly, if the Hessian matrix $\partial^w_{y_iy_j}h(x,y)$ must be in an interval $\left[\partial^2_{y_iy_j}h(x, \left[y\right])\right]$, and that interval only contains positive definite matrices, then $h(x, y)$ is concave up over the entire input interval $\left[y\right]$. The positive definite-ness of interval matrices can be evaluated by testing if the interval eigenvalues are all positive, or by Sylvester's criterion, a condition on the determinants of successively larger submatrices of the original.

\begin{theorem}
	\label{def:sylvester}
	\textbf{(Sylvester's Criterion)}
	
	A matrix $A\in\bbR^{n\times n}$ is positive-definite if and only if each of the leading principal minors (the upper left blocks of increasing size) of $A$ have positive determinants.
\end{theorem}
Unfortunately, computing the exact bounds of the interval determinant is an NP-hard problem   \cite{horacekDeterminantsIntervalMatrices2018}. The following branch-and-act algorithm avoids some of this cost by using larger intervals for the interval determinant, rather than its exact bounds.

\begin{figure}[h]
	\centering
	\includegraphics[width=0.45\textwidth]{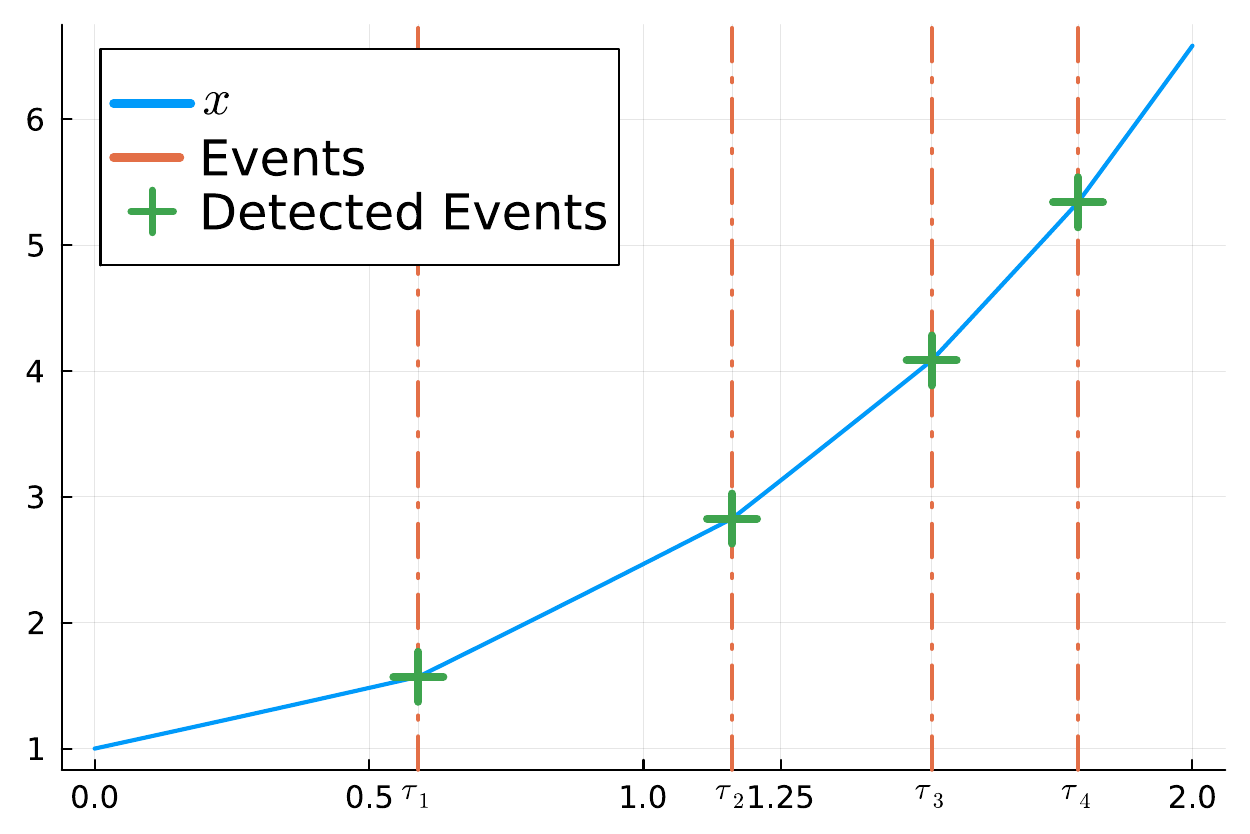}
	\includegraphics[width=0.45\textwidth]{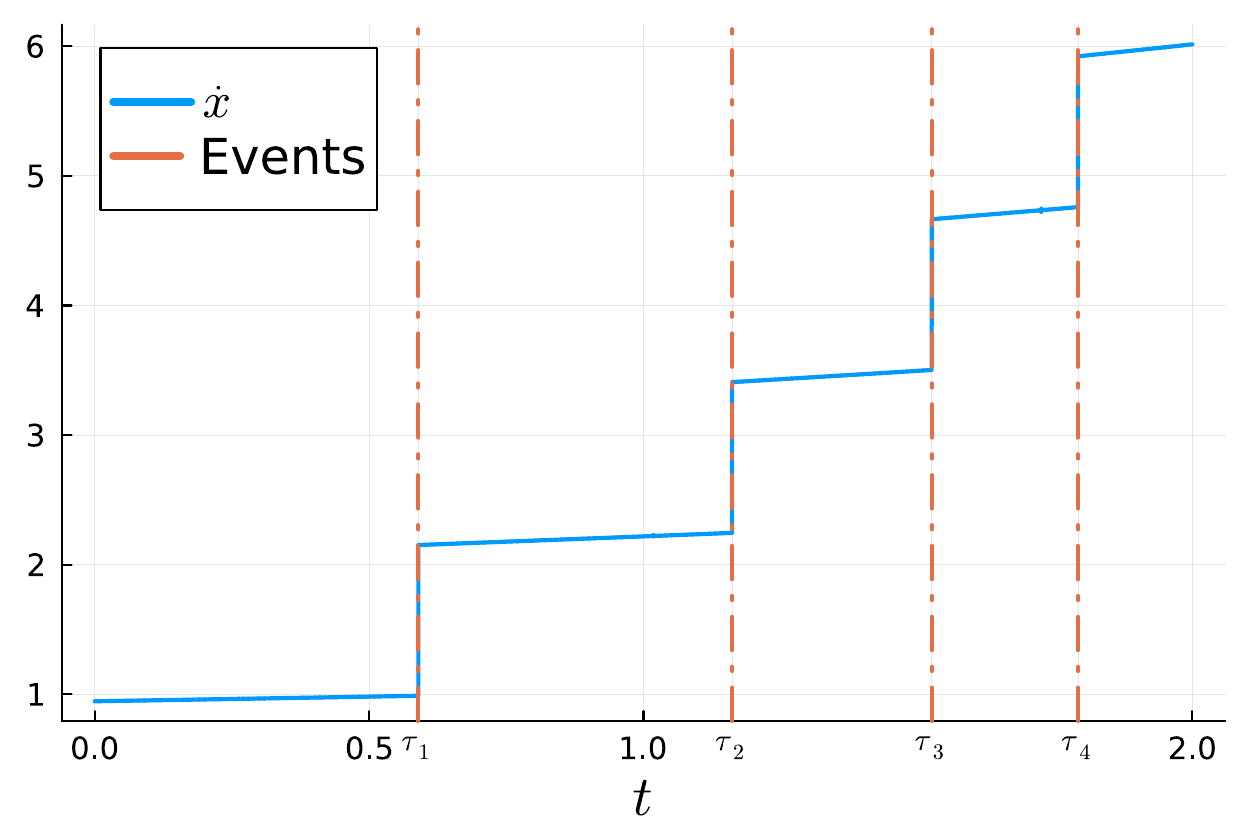}
	
	\caption{Results of running the solver on \eqref{eq:complicated-example}. The solver correctly detects the emergence and disappearance of local optimizers.}
	\label{fig:complicated-example-solution}
\end{figure}

\begin{algorithm}
	\label{algo:bnb-ia}
	\textbf{(Branch-and-Act with Interval Arithmetic)}
	
	Process the list of intervals to search for optimizers $\frS = \left\{[y]\right\}$ according to the following rules:
	\begin{enumerate} 
		\item Take the first item $[y]$ from $\frS$.
		\item \textbf{Evaluate Upper Bound:} Evaluate $[\ubar h, \bar h] = h(x, [y])$. If $\bar h$ is less than the current upper bound for $h$, update it. If $\ubar h$ is larger than the current upper bound for $h$, $[y]$ may be discarded if the desired behavior is to only find the global optimizer.
		\item \textbf{Gradient Test:} Evaluate the interval gradient $\partial_y h(x, [y])$. If the result interval does not contain $0$, $[y]$ contains no optimizers and can be discarded.
		\item \textbf{Hessian Test:} Test the interval Hessian $\partial^2_y h(x, [y])$ for positive definite-ness.
		\begin{enumerate}
			\item If the interval Hessian is negative definite, $h$ is concave down over the interval $[y]$, and $[y]$ can be discarded.
			\item If the interval Hessian is positive definite, $h$ is concave up over the interval $[y]$, and $[y]$ can be narrowed by any appropriate local optimization method.
		\end{enumerate}
		\item \textbf{Branch:} If the interval Hessian is neither positive- nor negative definite, decompose the interval $[y]$ and append the results to $\frS$. If the interval is bisected, this step can add at most $2^{n_y}$ new items to $\frS$.
		\item Repeat for all remaining items in $\frS$.
	\end{enumerate}
\end{algorithm}
This branch-and-act algorithm locates all intervals inside the initial search domain $\frS$ that contain minimizers of $h(x, y)$   \cite{rallGlobalOptimizationUsing1985, deussenGlobalSearchLocal2020}. The user is able to configure a termination criterion for the interval-narrowing and interval-splitting steps by setting the minimum width of an interval that will be split or narrowed. This directly controls the error in the optimizer's output. The solutions in \Cref{fig:easy-example-solution} and \Cref{fig:complicated-example-solution} were computed with a termination threshold of $\mathrm{width}\left[\ubar y^k, \bar y^k\right] \leq 10^{-8}$.

\section{Time Stepping}
The trapezoidal rule (a Runge-Kutta scheme of order 2   \cite{butcherNumericalMethodsOrdinary2008}) steps from time step $t^n$ to $t^{n+1}$ via an implicit formulation
\begin{equation}
	\label{eq:trap-rule}
	x^{n+1} = x^n + \frac{\Delta t}{2}(f(x^n, y^{\star,n})+ f(x^{n+1}, y^{\star,n+1}))
\end{equation}
Using this integrator requires estimates for $y^\star(t^{n})$ and $y^\star(t^{n+1})$.

\subsection{Local Optimizer Tracking}
If the time step $\Delta t$ is small enough, the integration step for $x$ can be augmented with $n_y \cdot \card\left(\left\{y^k\right\}\right)$ additional equations that force the values $\left\{y^k\right\}$ at time step $t^{n+1}$ to satisfy the first-order optimality conditions:
\begin{equation}
	\label{eq:first-order-opt-at-tn}
	\partial_{y^k_i}h\left(x^{n+1}, y^{k,n+1},_i\right) = 0
\end{equation}
	
\eqref{eq:trap-rule} and \eqref{eq:first-order-opt-at-tn} can be combined into one system of equations:
\begin{equation}
	\label{eq:integrator-with-tracking}
	\begin{aligned}
		0 &= x^n - x^{n+1} + \frac{\Delta t}{2}\left(\begin{aligned}&f(x^n, y^{\star,n})\,+\\&\qquad f(x^{n+1}, y^{\star,n+1})\end{aligned}\right)\\
		0 &= \partial_{y^k_i}h\left(x^{n+1}, y^{k,n+1},_i\right)
	\end{aligned}
\end{equation}
\eqref{eq:integrator-with-tracking} requires an initial guess for each $y^{k, n+1}$. With an estimate for $d_ty$, a suitable guess would be $y^{k, n+1} = y^{k, n} + d_ty^{k, n}\cdot\Delta t$. $d_ty$ can be computed from \eqref{eq:first-order-opt-at-tn} and the Implicit Function Theorem   \cite{deussenNumericalSimulationDifferentialalgebraic2023}:
\begin{equation*}
	\label{eq:local-tracking-guess}
	\begin{aligned}
	0&=\partial_{y}h(x, y^k)\\
	0&=d_x\partial_yh(x, y^k)\\
	0&=\partial^2_{yy}h(x, y^k)\cdot\partial_xy^k + \partial^2_{xy}h(x,y^k)\\
	\partial_xy^k &= -\left(\partial^2_{yy}h(x,y^k)\right)^{-1}\partial^2_{xy}h(x, y^k)
	\end{aligned}
\end{equation*}

Local optimizer tracking allows the solver to avoid expensive calls to global optimization while maintaining an estimate of each local optimizer. The error in this estimation is bounded by $\dot y^k(t)$, and it is corrected by periodically re-running the global optimizer. Re-running the global optimizer is also necessary to correctly identify newly-emerged or vanishing optimizers (see \eqref{eq:complicated-example}).  

Wthout an explicit treatment of time steps that bracket a jump event, the order of convergence of \textit{any} chosen integrator will be reduced to 1   \cite{mannshardtOnestepMethodsAny1978}.

\subsection{Jump Event Detection}

Fortunately, explicit treatment of the jump events can restore the order of convergence of the chosen integrator.
\begin{Definition}
	\label{def:event-fn}
	\textbf{(The Event Function)}
	Between any two local optimizers $y^1, y^2 \in \left\{y^k\right\}$, define
	\begin{equation*}
		H(x, y^1, y^2) = h(x, y^1) - h(x, y^2)
	\end{equation*}
	with the note that, at any $\tau_j$, an event occurs between $y^1$ and $y^2$ if $H\left(x(\tau_j), y^1, y^2\right) = 0$.
\end{Definition}
In order to locate an event, the solver must check the identity of the current global optimizer after \textit{each} time step. If it's changed, a root-finding procedure can be used on the event function to find $\tau_j$.

\section{Numerical Experiments}
\label{section:numerical-experiments}
The basic capabilities and peculiarities of the solver were tested using a simple example with a known analytic solution, and a more complex example designed to test the robustness of the event correction procedure. 

The solver relies on the templating mechanisms in \texttt{c++}. The Boost Interval library   \cite{melquiondBoostIntervalLibrary2022} provides an implementation for intervals of arbitrary real-number-like types.

Linear algebra operations were provided by Eigen  \cite{guennebaudEigenV32010}, and automatic differentiation capability was provided by \texttt{dco/c++}  \cite{leppkesDerivativeCodeOverloading2016}. For these experiments the solver (written in C++), was compiled using \texttt{-O3 -DNDEBUG} flags.
\subsection{Performance Testing}
An analytically-solvable DAEO is
\begin{equation}
	\label{eq:the-easy-one}
	\begin{aligned}
		x(0) &= 1\\
		\dot x(t) &= -(2+y^\star(t))x\\
		\left\{y^k(t)\right\} &= \argmin_y h(x,y)\\
		h(x, y) &= (1-y^2)^2 - (x-\frac{1}{2})\sin\left(\frac{\pi y}{2}\right) 
	\end{aligned}
\end{equation}
Here, the function $h(x, y)$ admits two local optimizers $\left\{-1, 1\right\}$, with 
\begin{equation*}
	y^\star = \begin{cases}
		-1 & x<\frac{1}{2}\\
		1 & x \geq \frac{1}{2}
	\end{cases}
\end{equation*}
The solution for the DAEO is then
\begin{equation}
	\label{eq:easy-daeo-solution}
	x(t) = \begin{cases}
		\exp\left(-3t\right) & t < -\frac{1}{3}\ln\frac{1}{2} \\
		\exp\left(-t + \frac{2}{3}\ln\frac{1}{2}\right) & t \geq -\frac{1}{3}\ln\frac{1}{2}
	\end{cases}
\end{equation}
The solver should detect exactly one event at $\tau_1 = -\frac{1}{3}\ln\frac{1}{2}$.

Since there are no constraints on either $x$ or $y$ in \eqref{eq:the-easy-one}, the KKT conditions for this example reduce to
\begin{equation*}
	\mathrm{minimize\ in\ } y:\,h\left(x(t), y\right)
\end{equation*}

A suitable local optimization algorithm will converge to one of the local optimizers. This result may be the global optimizer, but there is no way to know when the jump event occurs without a global optimizer. However, the branch-and-act algorithm described in \Cref{algo:bnb-ia} finds both local optimizers.

\begin{figure}[tb]
	\centering
	\includegraphics[width=0.45\textwidth]{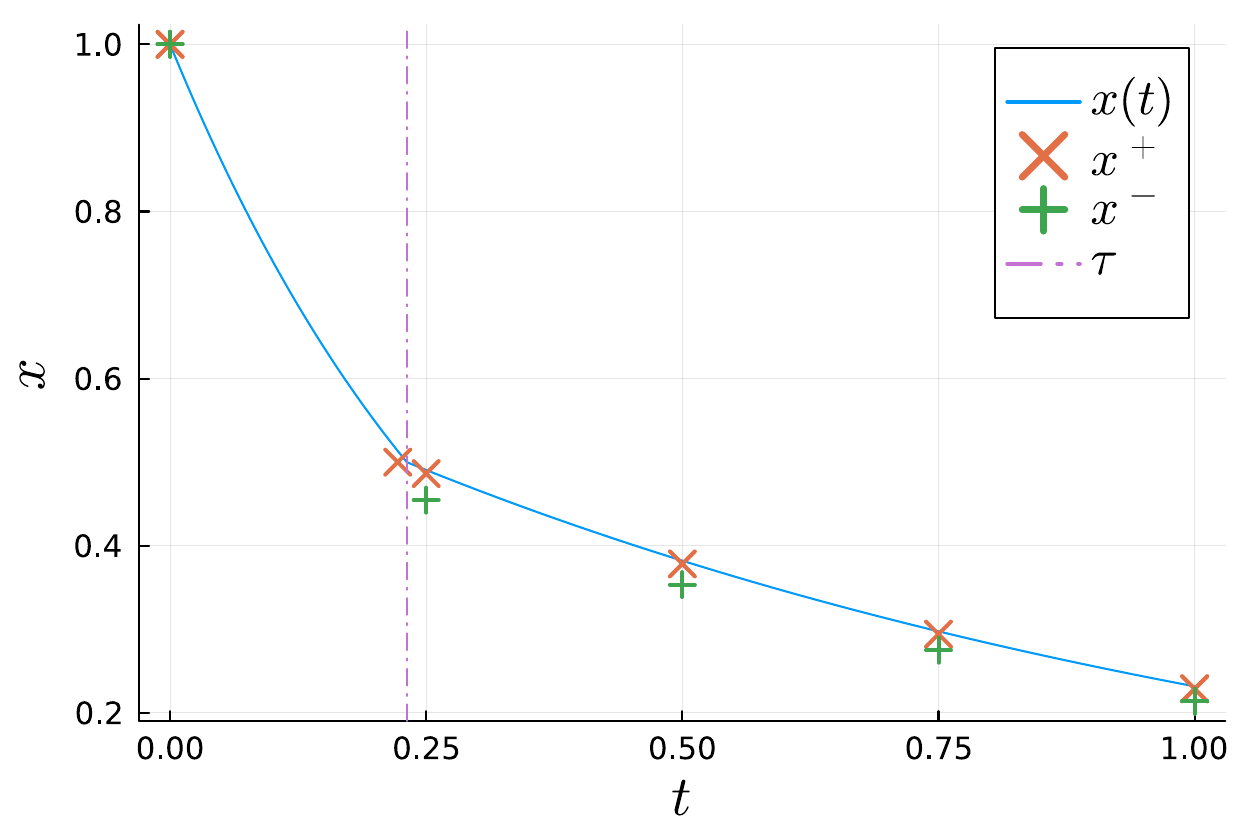}
	\caption{The solution to \eqref{eq:the-easy-one} as computed by the solver both with and without event correction enabled. $x(t)$ plots \eqref{eq:easy-daeo-solution}, and the vertical line is located at $\tau_1$, the location of the first event. The data points $x^+(t)$ and $x^-(t)$ show the computed trajectory for $x(t)$ with and without events, respectively.  The event detection procedure generates an extra data point at $t={\tau_1}$ that lies on the boundary between two of the component IVPs in \eqref{eq:the-easy-one}. $\Delta t$ was set to $0.25$ for this computation.}
	\label{fig:easy-example-solution}
\end{figure}
Figure \ref{fig:easy-example-solution} shows the difference between the solver results with the event detection and correction procedure enabled and without. The detected event adds an extra point into the computed solution, which functions as the initial value for a second DAE that exists from $t_{e} \leq t \leq t_{end}$. By correctly identifying the dynamics of $x(t)$ inside of the first time step from $t_0$ to $t_1$, the second-order convergence of the trapezoidal rule can be recovered.

\begin{figure}[h]
	\centering
	\includegraphics[width=0.45\textwidth]{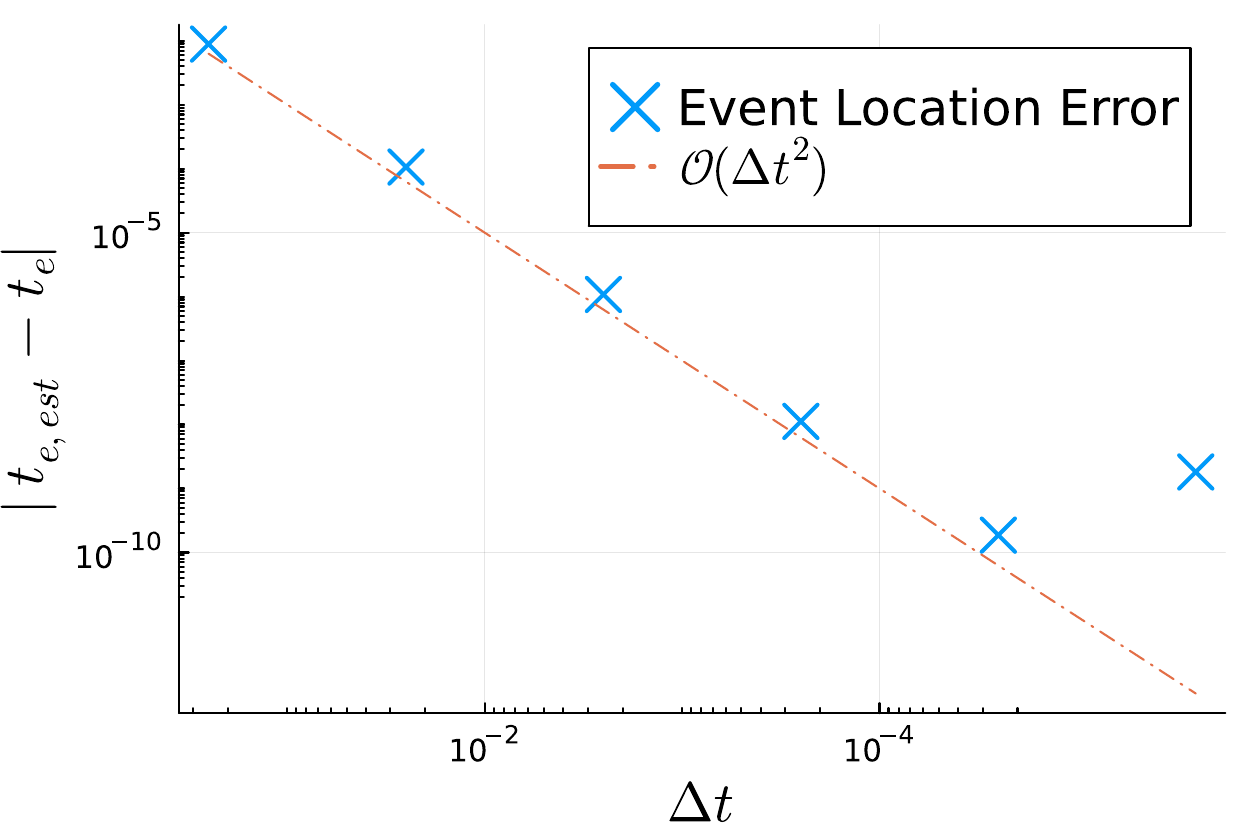}
	\caption{Event locator error $\abs{\tau - \tau_{exact}}$. The event locator error converges with order $\Delta t^2$ to a user-set tolerance, which is used for every root-finding routine in the solver.}
	\label{fig:easy-event-location}
\end{figure}
The event location procedure converges quickly enough to satisfy the consistency requirements given in   \cite{mannshardtOnestepMethodsAny1978}. Consequently, using the event detection and correction procedure will restore the trapezoidal rule to order 2. A user-specified parameter controls the termination criterion for the event-location procedure, which places a strict limit on its accuracy.

\begin{figure}[h]
	\centering
	\includegraphics[width=0.45\textwidth]{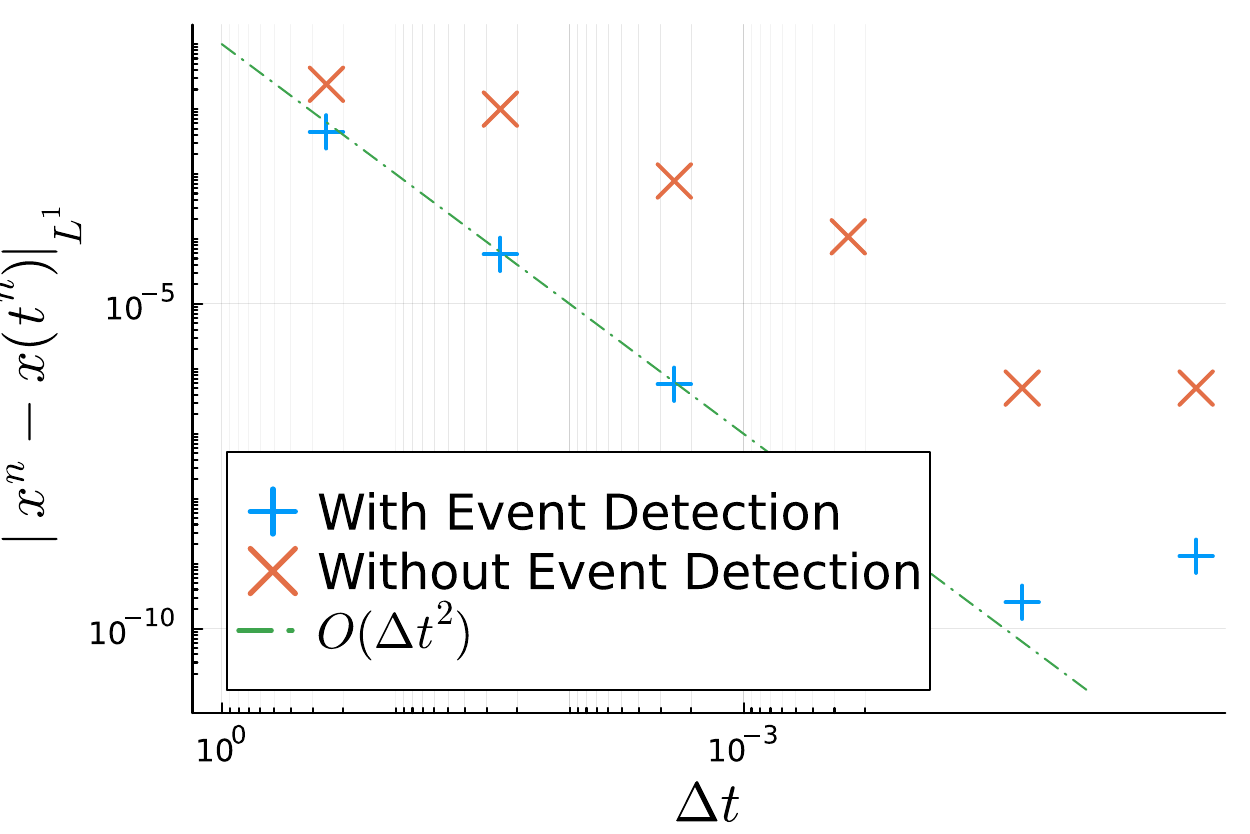}
	\caption{Convergence rates of the solver with (marked in blue crosses, converges with order 2) and without (marked in orange exes, converges with order 1).}
	\label{fig:easy-error-improvement}
\end{figure}
Figure \ref{fig:easy-error-improvement} shows the order-of-convergence increase due to the event location and simulation procedure. The price of recovering the second-order convergence is the addition of the event-detection procedure. The user may also trigger a search of the optimizer space at specific intervals, which incurs significant computational cost, especially if $n_y > 1$. However, local-optimizer tracking minimizes the number of times the global optimizer needs to be run. If it is known that the total number of optimizers $\left\{y^k\right\}$ is constant, the global optimizer only needs to be run as the beginning, which is the case for \eqref{eq:the-easy-one}.

\begin{figure}[h]
	\centering
	\begin{tabularx}{0.5\textwidth}{| >{\centering\arraybackslash}X | >{\raggedleft\arraybackslash}X |  >{\raggedleft\arraybackslash}X | >{\raggedleft\arraybackslash}X |}
		\hline
		$\Delta t$ & No Event Detection & Event Tracking & Always Optimize \\
		\hline
		2.50e-01 & 3 ms & 15 ms & 23 ms  \\
		2.50e-02 & 14 ms & 20 ms & 83 ms  \\
		2.50e-03 & 122 ms & 126 ms & 793 ms  \\
		2.50e-04 & 1125 ms & 1051 ms & 7593 ms  \\
		2.50e-05 & 11354 ms & 10444 ms & 74874 ms  \\
		2.50e-06 & 107298 ms & 109945 ms & 744696 ms \\\hline
	\end{tabularx}
	\caption{Run time of the solver (ms) compared against the time step size. The cost of global optimization at every time step, even in one dimension, quickly becomes prohibitive as $\Delta t$ decreases. There is little difference between the cost of the solver with event correction enabled or disabled, since the solution to \eqref{eq:the-easy-one} only involves one event correction.}
	\label{fig:easy-cost-comparison}
\end{figure}
Figure \ref{fig:easy-cost-comparison} shows the runtime costs of solving \eqref{eq:the-easy-one} in various configurations. When set to "Only Global Optimization" or "Always Optimize", the solver performs the optimization procedure (\ref{algo:bnb-ia}) in every time step and uses the results to perform event correction. The cost of this procedure scales significantly faster than the cost of the integrator alone. 
% explain why brancing works the way it does somewhere
In fact, the cost of the branch-and-act procedure scales as $2^{n_y}$ (each "branch" step generates $2^{n_y}$ new intervals to search), while solving the local optimizer only scales fast as the cost of solving a linear system of $n_x + n_y\cdot\card\left(\left\{y^k\right\}\right)$ equations. When configured as "No Event Detection", the solver skips the event location procedure, which only saves the cost of one root-finding operation and one integration step.

\subsection{Robustness Testing} The following example was designed to test the robustness of the solver in a situation where local optimizers emerge and vanish. 
\begin{equation}
	\label{eq:complicated-example}
	\begin{aligned}
		x(0) &= 1\\
		\dot x(t) &= y^\star(t)\\
		\left\{y^k(t)\right\} &= \argmin_{y} h(x, y)\\
		h(x, y) &= (x-y)^2 + \sin 5y
	\end{aligned}
\end{equation}
This particular objective function $h(x, y)$ has a set of local optimizers $\left\{y^k\right\}$ with a size that depends on the value of $x(t)$. 

The user is able to select the Lipshitz constants on $h$ that were derived in   \cite{deussenNumericalSimulationDifferentialalgebraic2023}. Correctly setting these coefficients requires some manual tuning, as incorrect values can lead to solutions that are subtly unreasonable (event detection in the wrong time step) or that are obviously incorrect (diverge to $\pm\infty$). Solutions to \eqref{eq:complicated-example} were computed with a minimum jump event size of $1.0\times10^{-4}$. Relaxing this setting any further generates nonsense results for $x(t)$, and reducing it causes the solver to erroneously detect events.

\section{Conclusions}
A new DAEO solver was developed to take advantage of the local optimizer tracking and jump event detection techniques presented in   \cite{deussenNumericalSimulationDifferentialalgebraic2023}. The solver was tested on example DAEO problems with both fixed- and time-dependent sets of optimizers $\left\{y^k\right\}$, and produces accurate results on both types of problems. However, the solver is very sensitive to parameters chosen by the user, namely, the decision threshold for event detection. The best possible choice for this threshold should never cause the erroneous detection of an event, but the possibility of emerging or vanishing optima makes this choice non-trivial.

Most significantly, the local optimizer tracking technique saves significant computational effort over current sequential methods for solving DAEOs  \cite{bieglerNonlinearProgrammingConcepts2010, plochDirectSingleShooting2022, plochMultiscaleDynamicModeling2019}, even if $n_y = 1$. 

\section{Further Work}

In order for the local optimizer tracking and event detection techniques to be relevant solution methods for larger problems ($n_x, n_y > 10$), such as the biorefinery simulation in \cite{plochMultiscaleDynamicModeling2019}, an alternate test for the positive definite-ness of the interval Hessian would have to be developed. Sylvester's Criterion (\ref{def:sylvester}) requires verifying the sign of the interval Hessian's determinant, which can be computed directly for smaller interval matrices. However, for larger interval matrices, directly computing the interval determinant in any reasonable amount of time requires decomposing the matrix, methods for which are developed in   \cite{goldsztejnGeneralizedIntervalLU2007, neumaierIntervalMethodsSystems1990a}. It would also be viable to take advantage of the symmetry of the interval Hessian, and directly test its interval eigenvalues   \cite{deifIntervalEigenvalueProblem1991}.

Currently, the Lipschitz coefficients and bounds for $\partial_ty^k$ are set by the user, and the solver is very sensitive to these chosen parameters. In some cases, it is feasible to analytically compute the maximum possible optimizer drift, and pass that as a specific parameter. The usefulness of the solver would, however, be improved if it were possible to estimate these parameters during run-time.

In the broader context of differentiable programming, the improvement of the order of convergence of DAEO solvers to match that of DAE solvers is an important step towards explicit analysis of the sensitivity of solutions to DAEOs to discontinuities in the solution of a global optimization problem. It is now more feasible to study the sensitivity of these solutions to input parameters and the positions of the jump events $t_{e_j}$, which may improve the accuracy and effectiveness of solutions to control problems, such as the UV-flash problem posed in   \cite{ritschelAlgorithmGradientbasedDynamic2018}.

\appendix
\section{Code Reference}
%% TODO upload and convert project to use AD.hpp
The code for this project is located on GitHub at
\href{https://github.com/STCE-at-RWTH/daeo-tracking-solver}{\texttt{STCE-at-RWTH/daeo-tracking-solver}}. 

It also relies on a public fork of Boost Interval that contains a patch developed by the authors, which is referenced via a \texttt{git submodule}.

% \newpage

\printbibliography
\end{document}